\documentclass[12pt]{amsart}

\usepackage[english]{babel}

\usepackage[utf8]{inputenc}

\usepackage{amsmath,nccmath}
\usepackage{amsfonts,amssymb,latexsym,graphics,amsthm}
\usepackage[bookmarks=false]{hyperref}
\usepackage[dvips]{color}
\usepackage[dvips]{graphicx}
\usepackage{pstricks,rotating}
\usepackage{url,multirow}

\renewenvironment{itemize}{\begin{list}{\labelitemi}{\leftmargin=1.5em}}{\end{list}}
\renewcommand{\labelitemi}{$\bullet$}

\newcommand{\ket}[1]{\ensuremath{|#1\rangle}}
\newcommand{\bra}[1]{\ensuremath{\langle #1|}}
\newcommand{\W}{\bra W}
\newcommand{\V}{\ket V}
\newcommand{\braket}[2]{\ensuremath{\langle #1|#2 \rangle}}

\newcommand{\x}{\bar x}
\newcommand{\y}{\bar y}
\newcommand{\z}{\bar z}

\newcommand{\hr}{{\rm hr}}
\newcommand{\hc}{{\rm hc}}
\newcommand{\fr}{{\rm fr}}
\newcommand{\fc}{{\rm fc}}
\newcommand{\emr}{{\rm emr}}
\newcommand{\emc}{{\rm emc}}
\newcommand{\fnr}{{\rm fnr}}
\newcommand{\fnc}{{\rm fnc}}
\newcommand{\dco}{{\rm dco}}
\newcommand{\lco}{{\rm lco}}
\newcommand{\mpx}{{\rm mp}}
\newcommand{\mi}{{\rm mi}}
\newcommand{\sd}{{\rm sd}}
\newcommand{\snd}{{\rm snd}}
\newcommand{\fd}{{\rm fd}}
\newcommand{\fnd}{{\rm fnd}}

\newcommand{\da}{\downarrow}
\newcommand{\la}{\leftarrow}

\newcommand{\bla}{
  \psset{unit=1mm} \begin{pspicture}(0,0) \rput(0,2){$\ddots$}   \end{pspicture} 
}

\title[Generalized Dumont-Foata polynomials and alternative tableaux]
  {Generalized Dumont-Foata polynomials and alternative tableaux}
\date{}
\author{Matthieu Josuat-Vergès}
\thanks{Partially supported by the grant ANR08-JCJC-0011.}
\address{LRI, Université Paris-Sud, Bât. 490, 91405 Orsay, FRANCE}
\email{josuat@lri.fr}

\hypersetup{
    unicode=false,          
    pdftoolbar=true,        
    pdfmenubar=true,        
    pdffitwindow=false,     
    pdfstartview={FitH},    
    pdftitle={Generalized Dumont-Foata polynomials and alternative tableaux}, 
    pdfauthor={Matthieu Josuat-Vergès (Université Paris XI)},
    pdfsubject={},
    pdfcreator={Creator},   
    pdfproducer={Producer}, 
    pdfkeywords={},
    pdfnewwindow=true,      
    colorlinks=true,       
    linkcolor=red,          
    citecolor=blue,        
    filecolor=magenta,      
    urlcolor=cyan           
}

\newtheorem{thm}{Theorem}[section]

\newtheorem{lem}[thm]{Lemma}
\newtheorem{prop}[thm]{Proposition}

\theoremstyle{definition}
\newtheorem{defn}[thm]{Definition}
\newtheorem{rem}[thm]{Remark}

\renewcommand{\arraystretch}{2}

\begin{document}

\begin{abstract}
Dumont and Foata introduced in 1976 a three-variable symmetric refinement of Genocchi 
numbers, which satisfies a simple recurrence relation. A six-variable generalization 
with many similar properties was later considered by Dumont. They generalize
a lot of known integer sequences, and their ordinary generating function can be 
expanded as a Jacobi continued fraction.

We give here a new combinatorial interpretation of the six-variable polynomials in terms 
of the alternative tableaux introduced by Viennot. A powerful tool to enumerate alternative 
tableaux is the so-called ``matrix Ansatz'', and using this we show that our combinatorial
interpretation naturally leads to a new proof of the continued fraction expansion.
\end{abstract}

\maketitle

\section{Introduction}

%

The unsigned Genocchi numbers $\{G_{2n}\}_{n\geq 1}$ can be defined through their generating 
function:
\begin{equation} 
  \sum_{n=1}^\infty G_{2n} \frac {x^{2n}}{(2n)!} = x\cdot \mathrm{tan}\left( \frac x2 \right).
\end{equation}
They are related with even Bernoulli numbers $B_{2n}$ by $G_{2n}=2(4^n-1)|B_{2n}|$, and they 
have a wide range of combinatorial properties \cite{DD95,AR94,Vie81,JZ96}. In the context of 
previous work by Carlitz, Riordan and Stein, an extension of these integers was proposed by 
Dumont and Foata \cite{DUFO}. It is defined by the recurrence $F_1(x,y,z)= 1$ and
\begin{equation} \label{df_rec}
  F_n(x,y,z)  =  (x+y)(x+z) F_{n-1}(x+1,y,z) -  x^2  F_{n-1}(x,y,z).
\end{equation}
They show that the polynomial $F_n$ is symmetric in $x$, $y$, and $z$, with non-negative 
coefficients, and such that $F_n(1,1,1)=G_{2n+2}$. Another nice property is that the generating 
function $\sum_{n=1}^{\infty} F_nt^n$ can be expanded as a J-fraction (more precision will be
given further in this introduction). Gessel and Zeng \cite{GZ} showed that $F_{n+1}$ is the $n$th 
moment of some orthogonal polynomials known as {\it continuous dual Hahn polynomials}, which are 
an important sequence in the Askey-Wilson hierarchy \cite{KoSw98}.

\bigskip

A further generalization of Genocchi numbers with many similar properties was defined by Dumont in 
terms of some combinatorial objects called {\it escaliers} \cite{DD95}. It is a sequence of 
six-variable polynomials $\Gamma_n(x,y,z,\x,\y,\z)$, or just $\Gamma_n$ for short. They can be 
characterized by a recurrence relation which generalizes \eqref{df_rec}, and has been obtained
independently by Randrianarivony \cite{AR94} and Zeng \cite{JZ96}. For brevity, let $\Gamma^+_n$ 
denote $\Gamma_{n}(x+1,y,z,\x+1,\y,\z)$.

\begin{defn} \label{def_dfg}
The generalized Dumont-Foata polynomials are defined by $\Gamma_1=1$ and
\begin{equation} \label{dfg_rec}
\Gamma_n  = (x+\z)(y+\x) \Gamma_{n-1}^+ + \big( x(\y-y) + \x(z-\z) - x\x \big) \Gamma_{n-1}.
\end{equation}
\end{defn}

Quite a lot of known integer sequences appear as specializations of $\Gamma_n$
\cite{AR94,RZ94,JZ96}:  Genocchi numbers, median Genocchi numbers, Euler numbers, 
median Euler numbers, Springer numbers. This polynomial $\Gamma_n$ generalizes
$F_n$ since we have $\Gamma_n(x,y,z,x,y,z)=F_n(x,y,z)$. 

Dumont \cite{DD95} conjectured that we have the following J-fraction for 
$\sum \Gamma_n t^n$:
\begin{equation} \label{dfg_frac}
  \sum_{n=1}^\infty \Gamma_n t^n = \cfrac{t}{1 - b_0t -  
        \cfrac{\lambda_1 t^2}{1 - b_1t -  \cfrac{\lambda_2t^2}{\ddots
}}},
\end{equation}
where the parameters $b_n$ and $\lambda_n$ are defined by:
\begin{equation}  \label{defbn}
\begin{split} 
b_n       & = (x+n)(\y+n)+(y+n)(\z+n)+(z+n)(\x+n) - n(n+1), \\
\lambda_n & = n(\x+y+n-1)(\y+z+n-1)(\z+x+n-1). 
\end{split}\end{equation}
This was proved independently by Randrianarivony \cite{AR94} and Zeng \cite{JZ96}
(and of course this implies the J-fraction expansion for $\sum F_nt^n$).
More precisely, Randrianarivony's method consists in the study of a {\it Stieltjes tableau}
and Zeng's method consists in calculations of {\it Hankel determinants}.

\bigskip

The main goal of this article is to give a new combinatorial interpretation of $\Gamma_n$ in 
terms of {\it alternative tableaux} \cite{PN09,Vie08} and six statistics on them, and obtain 
as a consequence a new proof of the continued fraction expansion \eqref{dfg_frac}. Alternative
tableaux were introduced by Viennot \cite{Vie08} in the context of a model of statistical physics
called Partially Asymmetric Simple Exclusion Process (PASEP), and previous work of Corteel
and Williams \cite{CW}. The ``matrix Ansatz'' first appeared in \cite{DEHP}, as a way to 
obtain the stationary distribution of the PASEP. In the combinatorial context, it is a method 
to enumerate these alternative tableaux in terms of operators satisfying certain relations.
We will also describe an analog of the matrix Ansatz to enumerate escaliers, which are the 
combinatorial objects used by Dumont to define $\Gamma_n$ in \cite{DD95}.

\bigskip

This article is organized as follows. In Section \ref{alt} we give definitions and 
known facts about alternative tableaux and the matrix Ansatz. In Section \ref{rec} we prove 
the new combinatorial interpretation of $\Gamma_n$ in terms of alternative 
tableaux using the recurrence \eqref{dfg_rec}. Section \ref{frac} contains our new proof of 
the continued fraction expansion \eqref{dfg_frac}. In Section \ref{sec_esc}, we describe 
the analog of the matrix Ansatz to enumerate escaliers.

\section{Alternative tableaux}
\label{alt}

Throughout this article we use the French convention for Young diagrams,
and Young diagrams may contain rows or columns of size 0. Any Young diagram is 
characterized by its upper-right boundary, which is a sequence of unit steps
going left or going down. We will encode this sequence by a word in the two 
letters $D$ and $E$, so that $D$ corresponds to the step $\rightarrow$ and $E$ 
corresponds to the step $\downarrow$. For example, $DDEDE$ is the Young diagram
with two rows of respective lengths 2 and 3.

\begin{figure}[h!tp] \center \psset{unit=4mm}
\begin{pspicture}(0,1)(5,5)
                 \psline(0,1)(0,5)
\psline(0,1)(4,1)\psline(1,1)(1,4)
\psline(0,2)(4,2)\psline(2,1)(2,3)
\psline(0,3)(3,3)\psline(3,1)(3,3)
\psline(0,4)(1,4)\psline(4,1)(4,2)
\rput(0.5,1.5){$\downarrow$}
\rput(2.5,1.5){$\downarrow$}
\rput(2.5,2.5){$\leftarrow$}
\end{pspicture}
\hspace{1.5cm}
\begin{pspicture}(0,1)(5,5)
                 \psline(0,1)(0,5)
\psline(0,1)(5,1)\psline(1,1)(1,5)
\psline(0,2)(4,2)\psline(2,1)(2,5)
\psline(0,3)(4,3)\psline(3,1)(3,4)
\psline(0,4)(3,4)\psline(4,1)(4,3)
\psline(0,5)(2,5)\psline(5,1)(5,1)
\rput(0.5,1.5){$\leftarrow$}
\rput(0.5,3.5){$\leftarrow$}
\rput(0.5,4.5){$\leftarrow$}
\rput(1.5,2.5){$\downarrow$}
\rput(3.5,2.5){$\downarrow$}
\end{pspicture}
\hspace{1.5cm}
\begin{pspicture}(0,1)(5,5)
                 \psline(0,1)(0,5)
\psline(0,1)(5,1)\psline(1,1)(1,5)
\psline(0,2)(5,2)\psline(2,1)(2,5)
\psline(0,3)(4,3)\psline(3,1)(3,5)
\psline(0,4)(4,4)\psline(4,1)(4,4)
\psline(0,5)(3,5)\psline(5,1)(5,2)
\rput(1.5,3.5){$\leftarrow$}
\rput(0.5,2.5){$\leftarrow$}
\rput(3.5,1.5){$\leftarrow$}
\rput(1.5,2.5){$\downarrow$}
\end{pspicture}
\hspace{1.5cm}
\begin{pspicture}(0,0)(5,5)
\psline(0,0)(5,0)\psline(0,0)(0,5)
\psline(0,1)(5,1)\psline(1,0)(1,5)
\psline(0,2)(4,2)\psline(2,0)(2,4)
\psline(0,3)(3,3)\psline(3,0)(3,3)
\psline(0,4)(2,4)\psline(4,0)(4,2)
\psline(0,5)(1,5)\psline(5,0)(5,1)
\rput(0.5,3.5){$\leftarrow$}
\rput(1.5,2.5){$\leftarrow$}
\rput(0.5,1.5){$\leftarrow$}
\rput(2.5,2.5){$\downarrow$}
\rput(3.5,1.5){$\downarrow$}
\end{pspicture}
\caption{\label{alt_ex} Examples of alternative tableaux. }
\end{figure}
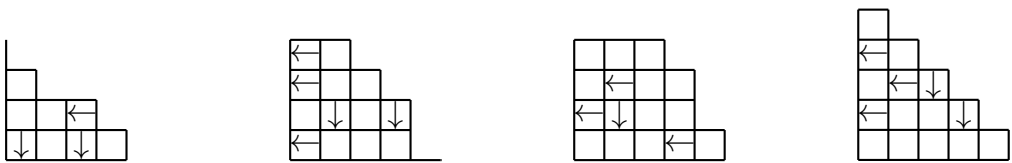

\begin{defn} 
Let $\lambda$ be a Young diagram. An {\it alternative tableau} of shape $\lambda$ is a filling 
of $\lambda$ such that each cell is either empty, contains an arrow $\la$ or an arrow $\da$, 
and each arrow has a clear view to the boundary. More precisely, all cells below a $\da$ in the 
same column (or to the left of a $\la$ in the same row) are empty. A column (respectively, 
row) of an alternative tableau is {\it free} if it contains no $\da$ (respectively, no $\la$). 
We denote by $\fr(T)$ (respectively, $\fc(T)$) the number of free rows (respectively, free 
columns) of $T$. See Figure~\ref{alt_ex} for examples. We use here the notation
with arrows as introduced by Nadeau \cite{PN09}.
\end{defn}

Alternative tableaux of a given shape can be enumerated via a method called matrix Ansatz. 
This method appeared in the context of a model of statistical physics (the partially asymmetric 
simple exclusion process), where it is used to derive the stationary probabilities of any state 
of the process.

\begin{prop}[Corteel-Williams \cite{CW}] \label{alt_ans}
Let $\bra{W}$ be a row vector, $\ket{V}$ a column vector, and $D$ and $E$ matrices such that:
\begin{equation} \label{de_rel}
\braket WV =1, \quad
\bra{W}E = \x \bra{W}, \qquad D\ket{V} = y \ket{V}, \quad \hbox{and} \quad DE-ED=D+E.
\end{equation}
Let $w$ be a word in the two letters $D$ and $E$, then we have
\begin{equation}
 \bra{W}w\ket{V} = \sum_{T} \x^{\fr(T)} y^{\fc(T)}
\end{equation}
where the sum is over alternative tableaux $T$ of shape $w$.
\end{prop}

The result of Corteel and Williams was actually stated in terms of {\it permutation tableaux},
which are slightly different objects. We refer to Viennot \cite{Vie08} and Nadeau \cite{PN09}
for the corresponding statement in terms of alternative tableaux, and for the bijection between
permutation tableaux and alternative tableaux. We have chosen to use alternative tableaux in 
this work because of their symmetry. Indeed, there is an elementary involution on alternative 
tableaux which is {\it conjugation}. To conjugate a tableau, take the image of the whole 
picture with respect to the South-West to North-East axis symmetry (in particular, the $\la$ 
and $\da$ are exchanged). See \cite{PN09} for details.

\medskip

Note that relations \eqref{de_rel} ensure that $\bra{W}w\ket{V}$ is well-defined and can
be computed explicitly. Indeed, we can use $DE-ED=D+E$ to obtain some coefficients
$c_{i,j}$ such that $w=\sum_{i,j} c_{i,j} E^iD^j$, and from the other relations 
we can obtain $\bra{W}w\ket{V}$. We refer to \cite{CW,CJW} for more details.

\medskip

Although not necessary to compute $\bra{W}w\ket{V}$ for a given word $w$, it is useful to 
have explicit matrices satisfying the PASEP matrix Ansatz. It can be checked \cite{DEHP} 
that the following $\mathbb{N}\times\mathbb{N}$-matrices:
\renewcommand{\arraystretch}{1.3}
\begin{equation} \label{def_DE}
D = \left(\begin{mmatrix}
y              &  1     &     &     &  (0)         \\
               &  y+1   & 2                        \\
               &        & y+2 & 3   & \phantom{a+2}\\ 
               &        &     & y+3 & \bla         \\
(0)               &        &     &     & \bla      \\
\end{mmatrix}\right)
,
E = \left(\begin{mmatrix}
\x    &         &         &        &  (0)    \\
y+\x  &  \x+1   \\
      &  y+\x+1 &  \x+2    &  \\ 
      &         &  y+\x+2 &  \x+3   &  \\
(0)   &         &         & \bla  & \bla \\
\end{mmatrix}\right),
\end{equation}
satisfy $DE-ED=D+E$. They are essentially a particular case of matrices defined by
Derrida \& al \cite{DEHP} in the context of the PASEP. As for the vectors, we can
take $\bra{W}=(1,0,0,\dots)$ and $\ket{V}=(1,0,0,\dots)^*$, and all relations in 
\eqref{de_rel} are satisfied.

\section{The new combinatorial interpretation of $\Gamma_n$}
\label{rec}

It is known \cite{Vie08} that $G_{2n+2}$ is the number of alternative tableaux whose shape
is the staircase with $n$ rows and columns, {\it i.e.} the Young diagram corresponding to 
the word $(DE)^n$. In \cite{CJW}, we have given three statistic in staircase alternative 
tableaux to give a combinatorial interpretation of $F_n(x,y,z)$. These are: the number 
of free rows, the number of free columns, and the number of corners containing $\la$ or
$\da$. Here, we give six statistics for the more general case of $\Gamma_n$.
Another difference is that in \cite{CJW}, the combinatorial interpretation was derived 
from the J-fraction for $\sum F_nt^n $, but here we use the recurrence relation
\eqref{dfg_rec} to prove the result.

\begin{defn} 
A column (respectively, row) of an alternative tableau is {\it empty} if it contains no $\da$ nor $\la$.
Let $T$ be an alternative tableau. We denote by:
\begin{itemize}
\item $\emr(T)$, the number of empty rows in $T$,
\item $\fnc(T)$, the number of free non-empty columns in $T$,
\item $\dco(T)$, the number of corners containing a $\da$ in $T$,
\item $\fnr(T)$, the number of free non-empty rows in $T$,
\item $\emc(T)$, the number of empty columns in $T$,
\item $\lco(T)$, the number of corners containing a $\la$ in $T$.
\end{itemize}
Moreover let $\mathcal{T}_n$ be the set of alternative tableaux whose shape is the staircase
Young diagram with $n$ rows and $n$ columns.
\end{defn}

For example, the rightmost tableau in Figure~\ref{alt_ex} is in $\mathcal{T}_5$, and the six
statistics that we have just defined are respectively 2, 2, 2, 0, 1, 0. The main new result of
this article is the following:

\begin{thm} \label{th_dfg}
For any $n\geq 1$, we have
\begin{equation} \label{dfg_alt}
  \Gamma_n(x,y,z,\x,\y,\z) = \sum_{T \in \mathcal{T}_{n-1}}
  x^{\emr(T)} y^{\fnc(T)} z^{\dco(T)} \x^{\fnr(T)} \y^{\emc(T)} \z^{\lco(T)}.
\end{equation}
\end{thm}

\begin{proof}
Both sides are equal to $1$ when $n=1$, so it suffices to show that the right-hand side 
satisfies the recurrence relation \eqref{dfg_rec}. We distinguish six kinds of tableaux 
in the set $\mathcal{T}_{n-1}$, according to the content of their leftmost column and 
upper left corner. Assuming that the theorem is true for $n-1$, we will show that these 
six kinds of tableaux have generating functions which add up to the right-hand side of
\eqref{dfg_rec}. This is summarized in the following table.

\smallskip

\begin{center}
\renewcommand{\arraystretch}{1.4}
\parbox{2.3cm}{ \vspace{5mm} The leftmost column is:} \quad
\begin{tabular}{c|c|c|c|}
           \multicolumn{4}{c}{ \qquad\qquad\qquad The upper left corner contains: }  \\[3mm]
         &   {$\da$}        &  {$\la$}    &  nothing    \\  
\hline
\multirow{2}{2.4cm}{ empty  }
&   \parbox{2cm}{\center $\times$ } & \parbox{2cm}{\center $\times$ } & Case 4 \\
                         &   &    &   $x\y\Gamma_{n-1} $     \\  
\hline
\multirow{2}{2.4cm}{ free non-empty }
&  \parbox{1.8cm}{\center $\times$ }  & Case 2 &  Case 5  \\
   &                 &   $y\z\Gamma_{n-1}^+$      &   $xy(\Gamma_{n-1}^+-\Gamma_{n-1})$ \\  
\hline
\multirow{2}{2.4cm}{ non-free }
  &  Case 1   &  Case  3     &  Case 6  \\
 & $\x z\Gamma_{n-1}$  &   $\x\z(\Gamma_{n-1}^+-\Gamma_{n-1})$ &  $x\x(\Gamma_{n-1}^+-\Gamma_{n-1})$ \\  
\hline
\end{tabular}
\end{center}

\smallskip

For example we will show that the tableaux of the fourth kind (case 4), {\it i.e.} those 
having an empty leftmost column, have generating function $x\y \Gamma_{n-1} $. The three 
cells containing a $\times$ in this table do not correspond to any tableaux.

\begin{figure}[h!tp]  \center \psset{unit=4mm}
\begin{pspicture}(0,0)(6,6)
\psline(0,0)(6,0)\psline(0,0)(0,6)
\psline(0,1)(6,1)\psline(1,0)(1,6)
\psline(0,2)(5,2)\psline(2,0)(2,5)
\psline(0,3)(4,3)\psline(3,0)(3,4)
\psline(0,4)(3,4)\psline(4,0)(4,3)
\psline(0,5)(2,5)\psline(5,0)(5,2)
\psline(0,6)(1,6)\psline(6,0)(6,1)
\psframe[linewidth=0.7mm](-0.1,-0.1)(1.1,6.1)
\rput(2.5,2.5){$\leftarrow$}
\rput(1.5,1.5){$\leftarrow$}
\rput(2.5,0.5){$\downarrow$}
\rput(4.5,1.5){$\downarrow$}
\rput(0.5,5.5){$\downarrow$}
\end{pspicture}
\hspace{1.2cm}
\begin{pspicture}(0,0)(6,6)
\psframe[fillstyle=solid,fillcolor=lightgray](0,0)(1,1)
\psframe[fillstyle=solid,fillcolor=lightgray](0,3)(1,5)
\psline(0,0)(6,0)\psline(0,0)(0,6)
\psline(0,1)(6,1)\psline(1,0)(1,6)
\psline(0,2)(5,2)\psline(2,0)(2,5)
\psline(0,3)(4,3)\psline(3,0)(3,4)
\psline(0,4)(3,4)\psline(4,0)(4,3)
\psline(0,5)(2,5)\psline(5,0)(5,2)
\psline(0,6)(1,6)\psline(6,0)(6,1)
\psframe[linewidth=0.7mm](-0.1,-0.1)(1.1,6.1)
\rput(2.5,2.5){$\leftarrow$}
\rput(1.5,1.5){$\leftarrow$}
\rput(2.5,0.5){$\downarrow$}
\rput(4.5,1.5){$\downarrow$}
\rput(0.5,5.5){$\leftarrow$}
\end{pspicture}
\hspace{1.2cm}
\begin{pspicture}(0,0)(6,6)
\psline(0,0)(6,0)\psline(0,0)(0,6)
\psline(0,1)(6,1)\psline(1,0)(1,6)
\psline(0,2)(5,2)\psline(2,0)(2,5)
\psline(0,3)(4,3)\psline(3,0)(3,4)
\psline(0,4)(3,4)\psline(4,0)(4,3)
\psline(0,5)(2,5)\psline(5,0)(5,2)
\psline(0,6)(1,6)\psline(6,0)(6,1)
\psframe[linewidth=0.7mm](-0.1,-0.1)(1.1,6.1)
\rput(2.5,2.5){$\leftarrow$}
\rput(1.5,1.5){$\leftarrow$}
\rput(2.5,0.5){$\downarrow$}
\rput(4.5,1.5){$\downarrow$}
\end{pspicture}
\hspace{1.2cm}
\begin{pspicture}(0,0)(6,6)
\psframe[fillstyle=solid,fillcolor=lightgray](0,0)(1,1)
\psframe[fillstyle=solid,fillcolor=lightgray](0,3)(1,5)
\psline(0,0)(6,0)\psline(0,0)(0,6)
\psline(0,1)(6,1)\psline(1,0)(1,6)
\psline(0,2)(5,2)\psline(2,0)(2,5)
\psline(0,3)(4,3)\psline(3,0)(3,4)
\psline(0,4)(3,4)\psline(4,0)(4,3)
\psline(0,5)(2,5)\psline(5,0)(5,2)
\psline(0,6)(1,6)\psline(6,0)(6,1)
\psframe[linewidth=0.7mm](-0.1,-0.1)(1.1,6.1)
\rput(2.5,2.5){$\leftarrow$}
\rput(1.5,1.5){$\leftarrow$}
\rput(2.5,0.5){$\downarrow$}
\rput(4.5,1.5){$\downarrow$}
\end{pspicture}
\caption{ \label{rectab}
Recursive construction of staircase alternative tableaux.
}
\end{figure}
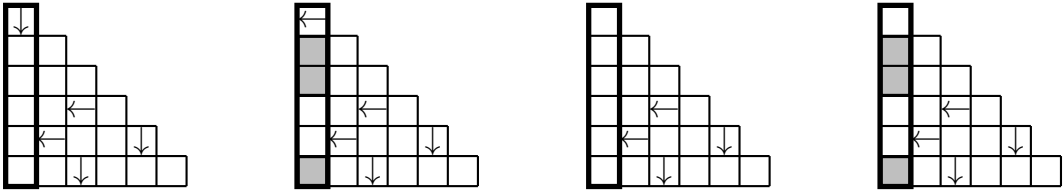

\begin{itemize}
\item Case 1. 
      When the upper-left corner contains a $\da$, there is no other arrow in the leftmost column.
      This corresponds to the first picture in Figure~\ref{rectab}. In this first kind of tableaux,
      the topmost row is free non-empty, and the upper left corner contains a $\da$, so this gives 
      a factor $\x z$. After removing the leftmost column, there can remain any tableau in 
      $\mathcal{T}_{n-2}$, hence the factor $\Gamma_{n-1}$. So the first kind of tableaux gives 
      indeed the term $\x z\Gamma_{n-1}$.
\item Case 2. 
      There is a factor $y\z$ since we assume that there is a $\la$ in the upper left corner, 
      and that the leftmost column is free and non-empty. These tableaux can be obtained the 
      following way: consider any tableau $T$ in $\mathcal{T}_{n-2} $, add a column to its left 
      with a $\la$ in the topmost cell of the added column. We color some of the other cells in 
      the added column in gray as in the second picture in Figure~\ref{rectab}, such that a cell
      is colored if there is no $\la$ to its right. Then, decide whether each gray cell is empty 
      or contains a $\la$. All tableaux of the second kind can be obtained this way, and the gray
      cells are in correspondence with the free rows of $T$. At the level of generating functions, this amounts to substitute $x$ with $x+1$ and $\x$ with $\x+1$. Indeed, an empty 
      (respectively, free non-empty) row remains so if we add nothing in the gray cell, but 
      becomes non-free if we add a $\la$.
\item Case 3. 
      This corresponds to the second picture in Figure~\ref{rectab}, but with the assumption that
      there is a $\da$ in one of the gray cells. Let us consider the set $S$ of tableaux of the 
      second kind such that there is at least a $\la$ in some gray cell. This set has generating 
      function $y\z(\Gamma^+_{n-1}-\Gamma_{n-1})$, indeed we have already $y\z\Gamma^+_{n-1}$ for 
      the all tableaux of the second kind and the term $-y\z\Gamma_{n-1}$ removes the cases where 
      all gray cells are empty. Then, there is a bijection between this set $S$ and the tableaux 
      of the third kind. Indeed, let $T\in S$, consider the bottommost $\la$ in the leftmost column
      of $T$, and replace this $\la$ with a $\da$. This way, we obtain exactly the tableaux of the 
      third kind. Replacing the $\la$ with $\da$ gives a factor $\x y^{-1}$ at the level of 
      generating functions. Thus we obtain $\x\z(\Gamma^+_{n-1}-\Gamma_{n-1})$ for the third kind 
      of tableaux.
\item Case 4. 
      This is similar to case 1, and corresponds to the third picture in Figure~\ref{rectab}.
      Here, removing the first column gives a factor $x\y$ since we assume that the leftmost 
      column is empty, and hence the upper row is empty. There can remain any tableau in 
      $\mathcal{T}_{n-2}$, hence the factor $\Gamma_{n-1}$.
      Thus we obtain $x\y\Gamma_{n-1}$ for the fourth kind of tableaux.
\item Case 5.
      This corresponds to the fourth picture in Figure~\ref{rectab}, with the assumption that
      the gray cells contain no $\da$ and at least a $\la$. Here the gray cells are obtained 
      exactly as in case 2 above. Proceeding similar to case 3 above, 
      we obtain $xy(\Gamma^+_{n-1}-\Gamma_{n-1})$ for the fifth kind of tableaux.
\item Case 6.
      There is a bijection between the fifth kind and the sixth kind of tableaux, similar to
      the bijection used in case 3. From a tableau of the fifth kind, consider the 
      bottommost $\la$ in the leftmost column, and replace it with a $\da$. Replacing the 
      $\la$ with $\da$ gives a factor $\x y^{-1}$ at the level of generating functions.
      Thus we obtain $x\x(\Gamma^+_{n-1}-\Gamma_{n-1})$ for the sixth kind of tableaux.
\end{itemize}
Adding the six terms in the above table, we get the right-hand side of \eqref{dfg_rec}. This 
shows that the right-hand side of \eqref{dfg_alt} satisfies the same recurrence as $\Gamma_n$, 
and completes the proof.
\end{proof}

As previously mentioned, there is a simple bijection between alternative tableaux and permutation
tableaux and it is possible to derive a combinatorial interpretation of $\Gamma_n$ in terms
of permutation tableaux, but the result is much more natural with the alternative tableaux.
In particular, the conjugation of alternative tableaux gives an easy way to prove a symmetry 
property of $\Gamma_n$, which has been first given by Randrianarivony \cite{AR94} and Zeng 
\cite{JZ96}. 
\begin{prop}[\cite{AR94,JZ96}] For any permutation $\sigma=u,v,w$  of $x,y,z$ we have:
\begin{itemize}
\item if $\sigma$ has signature 1, then $\Gamma_n(u,v,w,\bar u,\bar v,\bar w)
      =\Gamma_n(x,y,z,\x,\y,\z)$,
\item if $\sigma$ has signature -1, then $\Gamma_n(u,v,w,\bar u,\bar v,\bar w)
      =\Gamma_n(\x,\y,\z,x,y,z)$.
\end{itemize}
In particular, $F_n(x,y,z)=\Gamma_n(x,y,z,x,y,z)$ is symmetric in $x$, $y$, and $z$.
\end{prop}

\begin{proof}
From the recurrence relation \eqref{dfg_rec} , we have
\begin{equation} \label{sym_rec}
 \Gamma(x,y,z,\x,\y,\z) = \Gamma(\x,\z,\y,x,z,y).
\end{equation}
From the combinatorial interpretation in \eqref{dfg_alt} and using the conjugation of alternative 
tableaux, we have
\begin{equation} \label{sym_alt}
 \Gamma(x,y,z,\x,\y,\z) = \Gamma(\y,\x,\z,y,x,z).
\end{equation}
All symmetries given in the can be obtained by combining  \eqref{sym_rec} and
\eqref{sym_alt}.
\end{proof}

It is rather curious that one symmetry is obvious on the combinatorial interpretation, and
another one in the recurrence relation. In the model in terms of escaliers (see Section 
\ref{sec_esc}), only one symmetry is obvious, and it is both in the combinatorial interpretation 
and in the recurrence relation. Note that any symmetry of a generating function necessarily appear 
in the coefficients of its expansion as a J-fraction, and indeed it is straightforward to check 
that the coefficients $b_n$ and $\lambda_n$ defined in \eqref{defbn} have the same symmetries as 
$\Gamma_n$.

\section{The continued fraction expansion}
\label{frac}

In this section, we show that $\Gamma_n$ can be calculated via the matrix Ansatz, and we derive 
as a consequence a new proof of the continued fraction expansion for $\sum \Gamma_nt^n$. 
We consider the matrix
\begin{equation} \label{def_m}
M =  ED + (\z + x - \x)D + (z + \y - y)E + (\y-y)(x-\x)I,
\end{equation}
where $D$, $E$ are defined in \eqref{def_DE}, and $I$ is the identity matrix. It turns out that we can 
exploit Proposition~\ref{alt_ans} to obtain the following:

\begin{prop} \label{dfg_ans}
For any $n\geq0$, we have
\begin{equation}
\Gamma_{n+1}(x,y,z,\x,\y,\z) = \W M^n \V.
\end{equation}
\end{prop}

To prove this, we need a few helpful definitions and lemmas.

\begin{figure}[h!tp] \center \psset{unit=4mm}
\begin{pspicture}(0,0)(5,5)
\psline(0,0)(5,0)\psline(0,0)(0,5)
\psline(0,1)(5,1)\psline(1,0)(1,5)
\psline(0,2)(4,2)\psline(2,0)(2,4)
\psline(0,3)(3,3)\psline(3,0)(3,3)
\psline(0,4)(2,4)\psline(4,0)(4,2)
\psline(0,5)(1,5)\psline(5,0)(5,1)
\psline[linestyle=dashed,dash=0.8mm 0.8mm](1.5,0)(1.5,4)
\psline[linestyle=dashed,dash=0.8mm 0.8mm](0,0.5)(5,0.5)
\psline[linestyle=dashed,dash=0.8mm 0.8mm](4.5,0)(4.5,1)
\rput(0.5,3.5){$\leftarrow$}
\rput(0.5,1.5){$\leftarrow$}
\rput(2.5,1.5){$\downarrow$}
\rput(3.5,1.5){$\downarrow$}
\end{pspicture}
\hspace{1.5cm}
\begin{pspicture}(0,0)(5,5)
\psline(0,0)(5,0)\psline(0,0)(0,5)
\psline(0,1)(5,1)\psline(1,0)(1,5)
\psline(0,2)(4,2)\psline(2,0)(2,4)
\psline(0,3)(3,3)\psline(3,0)(3,3)
\psline(0,4)(2,4)\psline(4,0)(4,2)
\psline(0,5)(1,5)\psline(5,0)(5,1)
\psline[linestyle=dashed,dash=0.8mm 0.8mm](0.5,0)(0.5,5)
\psline[linestyle=dashed,dash=0.8mm 0.8mm](0,2.5)(3,2.5)
\rput(1.5,3.5){$\leftarrow$}
\rput(2.5,1.5){$\leftarrow$}
\rput(1.5,0.5){$\downarrow$}
\rput(4.5,0.5){$\downarrow$}
\end{pspicture}
\caption{\label{exttab} Example of two extended alternative tableaux.
}
\end{figure}
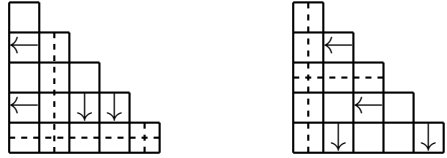

\begin{defn} Let $\mathcal{T}^*_n$ be the set of pairs $(T,X)$ where $T\in\mathcal{T}_n$ and $X$ 
is a subset of the empty rows and columns of $T$. Such a pair $(T,X)$ is called an {\it extended
tableau}, and will be represented the following way: from a picture of $T$, each row or column in 
$X$ is distinguished by a dashed line going through it. See Figure~\ref{exttab} for some examples.
Given $U=(T,X)\in\mathcal{T}^*_{n}$, we define
\begin{itemize}
\item $\hr(U)$ as the number of dashed rows,
\item $\hc(U)$ as the number of dashed columns.
\end{itemize}
For any statistic ``stat'' on alternative tableaux, and $U=(T,X)\in\mathcal{T}^*_n$ we define 
stat$(U)=$ stat$(T)$. For any extended tableau $U$, we define the weight $w(U)$ as 
\begin{equation}\begin{split}
 w(U) &= \x^{\emr(U)-\hr(U)} y^{\fnc(U)} z^{\dco(U)} \x^{\fnr(U)} y^{\emc(U)-\hc(U)} \z^{\lco(U)}
           (x-\x)^{\hr(U)} (\y-y)^{\hc(U)}  \\
      &= y^{\fc(U)-\hc(U)} z^{\dco(U)} \x^{\fr(U)-\hr(U)} \z^{\lco(U)}
           (x-\x)^{\hr(U)} (\y-y)^{\hc(U)},
\end{split}\end{equation}
the latest equality following from $\emr(T)+\fnr(T)=\fr(T)$ and $\emc(T)+\fnc(T)=\fc(T)$
for any alternative tableau $T$.
\end{defn}

\begin{lem} \label{geneprop}
We have:
\begin{equation} \label{gene}
 \sum_{U\in\mathcal{T}^*_n} w(U) = \sum_{T \in \mathcal{T}_{n}}
  x^{\emr(T)} y^{\fnc(T)} z^{\dco(T)} \x^{\fnr(T)} \y^{\emc(T)} \z^{\lco(T)}.
\end{equation}
\end{lem}
\begin{proof}
In the sum $\sum_{U\in\mathcal{T}^*_n} w(U)$, we have distinguished two
kinds of empty rows (dashed or non-dashed) with respective weights $x-\x$ and $\x$ instead
of one kind of empty row with weight $x$. Similarly we have distinguished two
kinds of empty columns (dashed or non-dashed) with respective weights $\y-y$ and $y$ instead
of one kind of empty column with weight $\y$. By an elementary argument, it is clear that these
distinctions do not change the generating function.
\end{proof}

\begin{defn}
Let $(T,X)$ be an extended tableau in $\mathcal{T}^*_n$. The {\it profile} of $(T,X)$ is the 
sequence $(i_1,\dots,i_n) $, where:
\begin{itemize}
\item $i_k=1$ if the $k$th corner of $T$ is empty,
\item $i_k=2$ if the $k$th corner of $T$ contains a $\la$,
\item $i_k=3$ if the $k$th corner of $T$ is in a dashed row but not in a dashed column, 
\item $i_k=4$ if the $k$th corner of $T$ contains a $\da$,
\item $i_k=5$ if the $k$th corner of $T$ is in a dashed column but not in a dashed row,
\item $i_k=6$ if the $k$th corner of $T$ is in a dashed column and in a dashed row.
\end{itemize}
Here the corners are numbered from the upper left one to the lower right one.
For example, the first extended tableau in Figure~\ref{exttab} has profile $(1,5,1,4,6)$,
and the second one has profile $(5,2,3,1,4)$.
\end{defn}

\begin{lem} \label{profile}
Let $M_1$, $\dots$, $M_6$ be the matrices
\begin{equation}
\begin{array}{lll} 
\label{def_mi}
M_1=ED,  \quad &M_2=\z D,    \quad &M_3=(x-\x)D,  \\
M_4=z E, \quad &M_5=(\y-y)E, \quad &M_6=(\y-y)(x-\x)I.
\end{array}
\end{equation}
For any $(i_1,\dots,i_n)\in\{1,\dots,6\}^n$, we have
\begin{equation} \label{profi}
 \sum_{U}  w(U) = \bra{W}M_{i_1}\dots M_{i_n} \ket{V},
\end{equation}
where the sum is over extended tableaux $U$ of profile $(i_1,\dots,i_n)$.
\end{lem}

\begin{proof}
Let $w$ be the word obtained from $i_1\dots i_n$ through the substitution $1\mapsto ED$, 
$2\mapsto D$, $3\mapsto D$, $4\mapsto E$, $5\mapsto E$, $6\mapsto \epsilon$ ($\epsilon$ being
the empty word). The main point is that there is a bijection $\phi$ between elements in 
$\mathcal{T}^*_n$ of profile $(i_1,\dots,i_n)$ and alternative tableaux of shape $w$. Indeed, 
to build an extended tableau, once the contents of the corners are specified it remains only 
to choose an alternative tableau of a smaller shape. More precisely, the bijection can be done 
the following way: 
\begin{itemize}
\item for each empty corner of the extended tableau, remove the corresponding cell in the 
      Young diagram,
\item shrink each dashed row or column, 
\item for each corner of the extended tableau containing a $\la$ (respectively, $\da$), 
      shrink the row (respectively, column) containing it.
\end{itemize}
See Figure~\ref{bij_phi} for an example, where we give the image of the two extended tableaux in
Figure~\ref{exttab}.

\begin{figure}[h!tp]  \center \psset{unit=4mm}
\hspace{1.5cm}
\begin{pspicture}(0,0)(3,4)
\psline(0,0)(2,0)\psline(0,0)(0,4)
\psline(0,1)(2,1)\psline(1,0)(1,3)
\psline(0,2)(1,2)\psline(2,0)(2,1)
\psline(0,3)(1,3)
\rput(0.5,2.5){$\leftarrow$}
\rput(0.5,0.5){$\leftarrow$}
\rput(1.5,0.5){$\downarrow$}
\end{pspicture}
\hspace{1.5cm}
\begin{pspicture}(0,0)(3,4)
\psline(0,0)(3,0)\psline(0,0)(0,3)
\psline(0,1)(3,1)\psline(1,0)(1,2)
\psline(0,2)(2,2)\psline(2,0)(2,2)
                 \psline(3,0)(3,1)
\rput(1.5,1.5){$\leftarrow$}
\rput(0.5,0.5){$\downarrow$}
\end{pspicture}
\caption{ \label{bij_phi}
Images of extended tableaux by the map $\phi$.}
\end{figure}
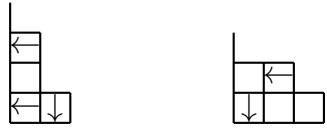

The weight of an extended tableau $U$ of profile $(i_1,\dots,i_n)$ is the product
\[
     y^a  z^b  \x^c  \z^d  (x-\x)^e  (\y-y)^f,
\]
where:
\begin{itemize}
\item $a$ is the number of non-dashed free columns in $U$, 
\item $b=\dco(U)$ is the number of 4's in $(i_1,\dots,i_n)$,
\item $c$ is the number of non-dashed free rows in $U$,
\item $d=\lco(U)$ is the number of 2's in $(i_1,\dots,i_n)$, 
\item $e=\hr(U)$ is the number of 3's plus the number of 6's in $(i_1,\dots,i_n)$,
\item $f=\hc(U)$ is the number of 5's plus the number of 6's in $(i_1,\dots,i_n)$.
\end{itemize}
An important property of the bijection $\phi$ is that the free rows (respectively, columns)
of $\phi(U)$ are in correspondence with non-dashed free rows (respectively, columns) of $U$.
It follows that
\begin{equation} \label{eq1}
    \sum_U w(U) =    z^b  \z^d  (x-\x)^e  (\y-y)^f \sum_T \x^{\fr(T)} y^{\fc(T)},
\end{equation}
where the first sum is over extended tableau of profile $(i_1,\dots,i_n)$ and the second
one is over alternative tableau of shape $w$.

Now, examine the product $M_{i_1}\dots M_{i_n}$. The factors $D$ and $E$ in this product readily 
gives the word $w$, and the other factors readily gives $ z^b \z^d (x-\x)^e  (\y-y)^f$, so:
\begin{equation} \label{eq2}
  \bra{W}M_{i_1}\dots M_{i_n} \ket{V} = z^b \z^d (x-\x)^e  (\y-y)^f \bra{W}w\ket{V}.
\end{equation}
Using Proposition~\ref{alt_ans}, the result follows from \eqref{eq1} and \eqref{eq2}.
\end{proof}

Now, we can prove Proposition~\ref{dfg_ans}.

\begin{proof}
Since $M=\sum_{i=1}^6 M_i$, the expansion of $M^n$ is also the sum of all products 
$M_{i_1}\dots M_{i_n}$ where $(i_1,\dots,i_n)$ runs through the set $\{1,\dots,6\}^n$. Hence,
\begin{equation}
 \bra{W}M^n\ket{V} = \sum_{(i_1,\dots,i_n)\in\{1,\dots,6\}^n} \bra{W} M_{i_1}\dots M_{i_n} \ket{V}.
\end{equation}
Using Equation~\eqref{profi} in Lemma~\ref{profile}, this gives
\begin{equation}
  \bra{W}M^n\ket{V} = \sum_{U\in\mathcal{T}^*_{n-1}} w(U).
\end{equation}
Using Equation~\eqref{gene} in Lemma~\ref{geneprop} and Theorem~\ref{th_dfg}, this is equal to 
$\Gamma_{n+1}$.
\end{proof}

From the definitions of $D$ and $E$ in \eqref{def_DE}, the matrix $M$ defined in \eqref{def_m}
can be calculated explicitly and we obtain the following statement.

\begin{prop} \label{calcm}
The matrix $M = (M_{i,j})_{i,j\in\mathbb{N}}$ is tridiagonal, and such that for any $i\geq0$ 
we have
\begin{equation}
M_{i,i} = b_i \quad \hbox{and} \quad M_{i,i+1}M_{i+1,i} = \lambda_{i+1},
\end{equation}
where $b_i$ and $\lambda_i$ are defined in \eqref{defbn}.
\end{prop}

\begin{proof}
We have
\begin{align*}
M_{i,i} &= E_{i,i}D_{i,i} + E_{i,i-1}D_{i-1,i} + (\z+x-\x) D_{i,i}+(z+\y-y) E_{i,i}+(\y-y)(x-\x) & \\
        &= (\x+i)(y+i) + (y+\x+i-1)i + (\z+x-\x)(y+i)  \\  & \qquad + (z+\y-y)(\x+i)  + (\y-y)(x-\x)   \\
        &= x\y + y\z + z\x + i \big( \x + \y + \z +x+y+z  \big) + i(2i-1) = b_i.
\end{align*}
We have also
\begin{align*}
M_{i,i+1} &= E_{i,i}D_{i,i+1} + (\z+x-\x) D_{i,i+1}   \\
          &= (\x+i)(i+1) + (\z+x-\x)(i+1) = (x+\z+i)(i+1),
\end{align*}
and
\begin{align*}
M_{i+1,i} &= E_{i+1,i}D_{i,i} + (z+\y-y) E_{i+1,i} = (y+\x+i)(y+i) + (z+\y-y)(y+\x+i) \\
          &= (\x+y+i)(z+\y+i).
\end{align*}
Hence, $M_{i,i+1}M_{i+1,i} = \lambda_{i+1}$.
It is straightforward to check that other coefficients in $M$ are 0, and this completes the proof.
\end{proof}

As a direct consequence of Propositions \ref{dfg_ans} and \ref{calcm}, let us give a new proof 
of the continued fraction expansion given in \eqref{dfg_frac}. First, note that $\bra{W}M^n\ket{V}$
is the top-left coefficient $(M^n)_{0,0}$ of the matrix $M^n$. This coefficient can be obtained
by expanding the product $M^n$ and we obtain
\begin{equation} \label{mat_dev}
\bra{W}M^n\ket{V} = \sum_{i_1,\dots,i_{n-1} \geq 0} M_{0,i_1}M_{i_1,i_2}\dots M_{i_{n-2},i_{n-1}} 
M_{i_{n-1},0}.
\end{equation}
Since the matrix $M$ is tridiagonal, we can restrict the sum to the case where two successive 
indices differ by at most 1, {\it i.e.} $|i_j-i_{j+1}|\leq1$ for any $j\in\{0,\dots,n-1\}$ where
$i_0=i_n=0$. These indices thus define the successive heights of a Motzkin path. Then \eqref{mat_dev}
shows that $\Gamma_{n+1}$ can be seen as the generating function of Motzkin paths of $n$ steps, with 
a weight $b_i$ for a level step at height $i$, and a weight $\lambda_i$ for a step $\nearrow$ between 
height $i-1$ and $i$. By a standard argument \cite{Fla} this implies the continued fraction 
given in \eqref{dfg_frac}.

\section{The matrix Ansatz for escaliers}

\label{sec_esc}

In the previous section, we have applied the link between alternative tableaux 
and matrices $D$ and $E$ satisfying $DE-ED=D+E$ to obtain the continued fraction.
In this section, we consider {\it escaliers}, which are the combinatorial objects
used by Dumont \cite{DD95} to define $\Gamma_n$. We will show that these objects 
can be enumerated by a similar method, but with matrices $B$ and $A$ satisfying 
$BA-AB=A+I$.

We will denote a Young diagram by a word in $B$ and $A$ in the same way that we
did with $D$ and $E$ ($B$ and $A$ respectively correspond to steps $\rightarrow$ 
and $\downarrow$ in the North-East boundary of the Young diagram).

\begin{figure}[h!tp] \center\psset{unit=4mm}
\begin{pspicture}(0,0)(8,4)
\psline(0,0)(8,0)  \psline(0,0)(0,4)
\psline(0,1)(8,1)  \psline(1,0)(1,4)
\psline(0,2)(6,2)  \psline(2,0)(2,4)
\psline(0,3)(4,3)  \psline(3,0)(3,3)
\psline(0,4)(2,4)  \psline(4,0)(4,3)
  \psline(5,0)(5,2)
  \psline(6,0)(6,2)
  \psline(7,0)(7,1)
  \psline(8,0)(8,1)
\rput(0.5,2.5){$\times$}\rput(1.5,3.5){$\times$}
\rput(3.5,2.5){$\times$}\rput(4.5,1.5){$\times$}\rput(2.5,0.5){$\times$}
\end{pspicture}
\hspace{2cm}
\begin{pspicture}(0,0)(8,4)
\psline(0,0)(8,0)  \psline(0,0)(0,4)
\psline(0,1)(8,1)  \psline(1,0)(1,4)
\psline(0,2)(6,2)  \psline(2,0)(2,4)
\psline(0,3)(4,3)  \psline(3,0)(3,3)
\psline(0,4)(2,4)  \psline(4,0)(4,3)
  \psline(5,0)(5,2)
  \psline(6,0)(6,2)
  \psline(7,0)(7,1)
  \psline(8,0)(8,1)
\rput(0.5,3.5){$\times$}
\rput(2.5,2.5){$\times$}
\rput(3.5,1.5){$\times$}\rput(4.5,0.5){$\times$}\rput(5.5,1.5){$\times$}
\rput(6.5,0.5){$\times$}
\end{pspicture}
\caption{ \label{pretab}
Examples of surjective pretableaux of shape $(BBA)^4$.
}
\end{figure}
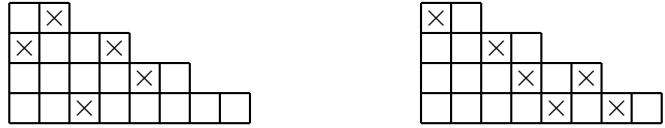

\begin{defn}
A {\it surjective pretableau} is a partial filling a Young diagram with $\times$, 
such that there is at least one $\times$ in each row and at most one $\times$ in each
column. A {\it surjective tableau} is a surjective pretableau such that there is exactly
one $\times$ in each column. An {\it escalier} (of size $n$) is a surjective tableau of 
shape $(BBA)^n$. See Figure~\ref{pretab} for some examples.
\end{defn}

The definition of $\Gamma_n$ by Dumont \cite{DD95} is given in terms of escaliers of size 
$n$ and six statistics on them. There is an obvious bijection between escaliers of size 
$n$ and surjective pretableau of shape $(BBA)^{n-1}$ (remove the bottom row of the escalier), 
so that Dumont's definition is equivalent to the following (as mentioned in the introduction, 
this is also known to be equivalent with Definition \ref{def_dfg}).

\begin{defn} A {\it co-corner} of the Young diagram $(BBA)^n$ is a cell which is
the left neighbor of a corner (for example the upper-left cell is a co-corner).
Let $\mathcal{S}_n$ be the set of surjective pretableau of shape $(BBA)^n$.  
A column is {\it empty} if it contains no $\times$. A $\times$ is {\it doubled}
if there is another $\times$ in the same row. We denote by:
\begin{itemize}
\item $\mi(T)$, the number of empty columns of odd index,
\item $\fd(T)$, the number of corners containing a doubled $\times$,
\item $\snd(T)$, the number of co-corners containing a non-doubled $\times$,
\item $\mpx(T)$, the number of empty columns of even index,
\item $\fnd(T)$, the number of corners containing a non-doubled $\times$,
\item $\sd(T)$, the number of co-corners containing a doubled $\times$,
\end{itemize}
Eventually, $\Gamma_n$ can be defined as 
\begin{equation}
  \Gamma_{n} = \sum_{T\in\mathcal{S}_{n-1}}
  x^{\mi(T)} y^{\fd(T)} z^{\snd(T)} \x^{\mpx(T)} \y^{\fnd(T)} \z^{\sd(T)}.
\end{equation}
\end{defn}

For example, the values of the six statistics on the first surjective pretableau in
Figure~\ref{pretab} are 1, 1, 1, 2, 1, and 0. As for the second one, the values are
0, 1, 2, 2, 0, and 1. The fact that these objects also follow the recurrence \eqref{dfg_rec} 
is seen by distinguishing several kinds of elements in $\mathcal{S}_{n-1}$ according to the 
content of the bottom row \cite{AR94,JZ96}.

The analog of the matrix Ansatz for escaliers is given in the following proposition.
The proof is similar to the case of alternative tableaux, and various examples of this kind 
of results were given in \cite{CJW}.

\begin{prop} \label{esc_an}
Let $\bra{W}$ be a row vector, $\ket{V}$ a column vector, and $A$ and $B$ matrices such that:
\begin{equation}  
\braket WV =1, \quad \bra{W}A = 0 , \quad  B\ket{V} = 0   , \quad BA-AB=A+I.
\end{equation}
Let $w$ be a word in the two letters $B$ and $A$, then the number of surjective tableaux of 
shape $w$ is $\bra{W}w\ket{V}$.
\end{prop}

\begin{proof}
This is done by a recurrence on the number of cells in the Young diagram. If there is no 
cell, then $w=A^iB^j$ for some $i$ and $j$, so $\bra{W}w\ket{V}$ equals 0 if $i>0$ or $j>0$ 
and 1 otherwise. Since there is at least a $\times$ in each row and column of a surjective 
tableau, there is no such tableau of shape $A^iB^j$ if $i>0$ or $j>0$, but we do have the 
``empty'' surjective tableau in the ``empty Young diagram'' when $i=j=0$.

Next, consider a word $w$ which is not in the form $A^iB^j$. It means we can factorize it 
in $w=w_1BAw_2$, and the factor $BA$ corresponds to a corner of the Young diagram. We can
distinguish three kinds of surjective tableaux of shape $w$, depending on the content of this
corner.
\begin{itemize}
\item If the corner is empty, we can remove it and obtain any surjective tableaux of shape 
      $w_1ABw_2$. By the recurrence assumption, their number is $\bra{W}w_1ABw_2\ket{V}$.
\item If the corner contains a doubled $\times$, we can delete the corner and its column, 
      and obtain any surjective tableau of shape $w_1Aw_2$. Their number is 
      $\bra{W}w_1Aw_2\ket{V}$.
\item If the corner contains a non-doubled $\times$, we can remove the corner and its row 
      and column, and obtain any surjective tableau of shape $w_1w_2$. Their number is 
      $\bra{W}w_1w_2\ket{V}$.
\end{itemize}
It follows that the number of surjective tableaux of shape $w$ is 
\[
  \bra{W}w_1ABw_2\ket{V} + \bra{W}w_1Aw_2\ket{V} + \bra{W}w_1w_2\ket{V} = \bra{W}w_1BAw_2\ket{V}
  = \bra{W}w\ket{V},
\]
and this completes the recurrence.
\end{proof}

It is easily checked that the following $\mathbb{N}\times\mathbb{N}$-matrices:
\renewcommand{\arraystretch}{1.2}
\begin{equation} \label{def_BA}
B = \left(\begin{mmatrix}
0    & 1 &   &   & (0)  \\
     & 1 & 2            \\
     &   & 2 & 3 &      \\ 
     &   &   & 3 & \bla \\
(0)  &   &   &   & \bla \\
\end{mmatrix}\right)
,\qquad
A = \left(\begin{mmatrix}
0   &   &   &   &  (0) \\
1   & 0                \\
    & 1 & 0            \\ 
    &   & 1 & 0 & \\
(0) &   &   & \bla  & \bla \\
\end{mmatrix}\right),
\end{equation}
satisfy $BA-AB=A+I$. We keep the definition of $\bra{W}$ and $\ket{V}$ as in the previous 
sections, since this also ensures that we have $ \bra{W}A = 0$ and  $B\ket{V} = 0 $.

To see how to use this result in the case of $\Gamma_n$, let us begin with the particular 
case $y=z=\y=\z=1$. We know that the number of surjective tableaux of shape $(BBA)^n$ is 
$\bra{W}(BBA)^n\ket{V}$. If we want to count surjective pretableaux, we have to authorize 
empty columns, and this is done by replacing $B$ with $B+I$. Indeed, in the expansion of 
the product $\bra{W}((B+I)(B+I)A)^n\ket{V}$, the choice of $B$ or $I$ in some factor 
corresponds to the choice of leaving a column empty or not. So the number of surjective 
pretableaux of shape $(BBA)^n$ is $\bra{W}((B+I)(B+I)A)^n\ket{V}$. If we want to follow 
the empty columns of odd (respectively, even) index by the parameter $x$ (respectively 
$\x$), it suffices to mark the terms $I$, and we obtain 
$\Gamma_{n+1}(x,1,1,\x,1,1)=\bra{W}((B+xI)(B+\x I)A)^n\ket{V}$.

As for the general case, in the same way that we have obtained Proposition~\ref{dfg_ans}
from the combinatorial interpretation in terms of alternative tableaux, we can obtain the 
following from the combinatorial interpretation in terms of escaliers.

\begin{prop} \label{dfg_ans2}
For any $n\geq0$, we have $\Gamma_{n+1} = \bra{W} N^n \ket{V}$ where $N$ is the matrix
\begin{equation}  \label{defN}
  N = A(B+xI)(B+\x I) + y\z(A+I) + (zI+\z A)(B+\x I) + (\y I +yA) (B+xI).
\end{equation}
\end{prop}

As in the previous section, we need some helpful definitions and lemmas.

\begin{defn}
The {\it profile} of $T\in\mathcal{S}_n$ is the sequence $(i_1,\dots,i_n)$, where
\begin{itemize}
\item $i_k=1$ if the $k$th co-corner and $k$th corner are empty,
\item $i_k=2$ if the $k$th co-corner and $k$th corner both contain a $\times$,
\item $i_k=3$ if the $k$th co-corner contains a non-doubled $\times$,
\item $i_k=4$ if the $k$th co-corner contains a doubled $\times$ and the $k$th corner is empty,
\item $i_k=5$ if the $k$th corner contains a non-doubled $\times$,
\item $i_k=6$ if the $k$th corner contains a doubled $\times$ and the $k$th co-corner is empty.
\end{itemize}
\end{defn}

For example, the two surjective pretableaux in Figure~\ref{pretab} have respective profile
$(5,6,3,1)$ and $(3,3,6,4)$.

\begin{lem} \label{profile2}
We define the matrices $N_1=A(B+xI)(B+\x I)$, $N_2=y\z(A+I)$, $N_3=z(B+\x I)$, 
$N_4=\z A (B+\x I)$, $N_5=\y (B+xI)$, $N_6=yA(B+xI)$.
For any $(i_1,\dots,i_n)\in\{1,\dots,6\}^n$, we have
\begin{equation} \label{profi2}
 \sum_{T} x^{\mi(T)} y^{\fd(T)} z^{\snd(T)} \x^{\mpx(T)} \y^{\fnd(T)} \z^{\sd(T)}
 = \bra{W}N_{i_1}\dots N_{i_n} \ket{V},
\end{equation}
where the sum is over $T\in\mathcal{S}_n$ of profile $(i_1,\dots,i_n)$.
\end{lem}

\begin{proof}
We follow the same scheme as in Proposition~\ref{profile}, and use Proposition~\ref{esc_an}.
Here, the surjective pretableaux of a given profile are in bijection with surjective
tableaux of a particular shape. Rather than giving a formal detailed proof, we sketch how
to understand the matrices $N_1$ through $N_6$, having in mind the proof of Proposition
\ref{esc_an} and the way surjective tableaux are build recursively. 

\begin{itemize}
\item The $k$th co-corner and $k$th corner correspond to the $k$th factor $BBA$ in 
      the word $(BBA)^n$. If these are empty ({\it i.e.} $i_k=1$), we can remove the 
      two cells and replace the factor $BBA$ with $ABB$. With the terms $xI$ and $\x I$ as
      seen before, we see that this case correspond to the matrix $N_1=A(B+xI)(B+\x I)$.
\item If $i_k=2$, {\it i.e.} the $k$th co-corner and $k$th corner both contain a $\times$, 
      we can remove the two columns containing these $\times$, and remove the factor $BB$ in 
      $BBA$. But we need to distinguish two cases, depending on wether there is a third 
      $\times$ in the same row or not, and if there is not, we also remove the row. 
      This gives a factor $A+I$, and there is a weight $y\z$ because of the doubled $\times$
      in the corner and co-corner. Hence this case correspond to the matrix $N_2=y\z(A+I)$.
\item If $i_k=3$, {\it i.e.} the $k$th co-corner contains a non-doubled $\times$, we can 
      remove its row and column, so the $k$th factor $BBA$ becomes a $B+\x I $. There is a 
      weight $z$ for the non-doubled $\times$ in the co-corner. Hence this case gives the 
      matrix $N_3=z(B+\x I)$.
\item If $i_k=4$, the difference with the previous case is that the $\times$ in the $k$th 
      co-corner is doubled, so we do not remove its row. So there remains a factor $A$, and
      there is a weight $\z$ instead of $z$. Hence this case gives the matrix 
      $N_4=\z A (B+\x I)$.
\item If $i_k=5$, this is similar to the case when $i_k=3$. But the weight is $\y$ instead
      of $z$ for the non-doubled $\times$ in the corner, and there remains a column of odd
      index so this gives a factor $B+xI$ instead of $B+\x I$. Hence this case gives the
      matrix $N_5=\y (B+xI)$.
\item If $i_k=6$, this is similar to the case when $i_k=4$. But the weight is $y$ instead
      of $\z$ for the doubled $\times$ in the corner, and there remains a column of odd
      index so this gives a factor $B+xI$ instead of $B+\x I$. Hence this case gives the
      matrix $N_6= yA(B+xI)$.
\end{itemize}
When we form the product $N_{i_1}\dots N_{i_n}$, it is clear that the matrix $N_{i_k}$ will
impose the conditions on the $k$th co-corner and corner, and hence 
$\bra{W}N_{i_1}\dots N_{i_n}\ket{V}$ is the generating function for elements in 
$\mathcal{S}_n$ of profile $(i_1,\dots,i_n)$.
\end{proof}

Now, we can prove Proposition~\ref{dfg_ans2}.

\begin{proof}
We have $N=\sum_{i=1}^6 N_i$, hence using Lemma~\ref{profile2}:
\begin{align*}
 \bra{W}N^n\ket{V} &= \sum_{(i_1,\dots,i_n)\in\{1,\dots,6\}^n} \bra{W}N_{i_1}\dots N_{i_n}\ket{V}\\
      &=  \sum_{T\in\mathcal{S}_{n}}
          x^{\mi(T)} y^{\fd(T)} z^{\snd(T)} \x^{\mpx(T)} \y^{\fnd(T)} \z^{\sd(T)} = \Gamma_{n+1}.
\end{align*}
This completes the proof.
\end{proof}

From the definitions of $B$ and $A$ in \eqref{def_BA}, the matrix $N$ defined in \eqref{defN}
can be calculated explicitly and we obtain the following statement.

\begin{prop}
The matrix $N = (N_{i,j})_{i,j\in\mathbb{N}}$ is tridiagonal, and such that for any $i\geq0$ 
we have $N_{i,i}=b_i$ and $N_{i,i+1}N_{i+1,i} = \lambda_{i+1}$.
\end{prop}

\begin{proof}
Straightforward calculations show that $N_{i,j}=0$ if $|i-j|>1$, $N_{i,i}=b_i$, and:
\begin{equation}
  N_{i,i+1}=(i+1)(z+\y+i), \qquad N_{i+1,i} = (x+\z+i)(y+\x+i).
\end{equation}
This gives indeed $N_{i,i+1}N_{i+1,i} = \lambda_{i+1}$.
\end{proof}

As in the case of alternative tableaux in the previous section, the previous two propositions means 
that the continued fraction expansion given in \eqref{dfg_frac} can be derived from the combinatorial 
interpretation in terms of escaliers, and the matrix Ansatz for escaliers given in
Proposition~\ref{esc_an}.

\begin{rem} Observe that $N\neq M$ since the non-diagonal coefficients are not the same,
but the two matrices are equal after a permutation of the variables $(x,y,z,\x,\y,\z)$.
However, there is {\it a priori} no simple way to link the matrices $M_i$'s with the $N_i$'s,
so we tend to think that despite the similarities in the method, there are two really 
different ways to obtain the continued fraction from a combinatorial model using the 
matrix Ansatz approach.

It is natural to ask if there is a bijection between $\mathcal{T}_{n}$ and $\mathcal{S}_{n}$
preserving the six statistics for $\Gamma_n$. From the fact that the recurrence relation is
checked on both sets, in theory it might be possible to describe recursively such a bijection.
It would be quite interesting to give a better answer to this question by providing a direct
bijection between our alternative tableaux and Dumont's escaliers.
\end{rem}

\section*{Conclusion}

Our main result is the new combinatorial interpretation of $\Gamma_n$ with the alternative
tableaux. We obtain two new proofs of Dumont's conjecture, with the same method applied to 
the two different combinatorial interpretations of $\Gamma_n$. What it is interesting
about these proofs is that they fit in the general framework settled in \cite{CJW},
linking J-fractions, operators satisfying certain commutation relations, and combinatorial 
objects.

\section*{Acknowledgement}

I thank Dominique Dumont, Dominique Foata, Jiang Zeng, Xavier Viennot and Sylvie Corteel 
for their conversation and encouragements.


\bigskip

\end{document}